\documentclass[11pt]{article}
\usepackage[a4paper,margin=25mm]{geometry}
\usepackage{amsmath,amssymb,mathtools,bm}
\usepackage[T1]{fontenc}
\usepackage{lmodern}
\usepackage{microtype}
\usepackage{enumitem}
\usepackage{booktabs}
\usepackage{array}
\usepackage{hyperref}
\usepackage{cleveref}
\usepackage{mathrsfs}
\allowdisplaybreaks
\hypersetup{colorlinks=true,linkcolor=black,citecolor=black,urlcolor=black,pdftitle={A Power Geometry and Hypergeometric Functions},pdfauthor={Shintaro Yoshizawa}}
\newtheorem{theorem}{Theorem}[section]
\newtheorem{proposition}[theorem]{Proposition}
\newtheorem{fact}[theorem]{Fact}
\newtheorem{definition}[theorem]{Definition}
\newtheorem{remark}[theorem]{Remark}
\newtheorem{problem}[theorem]{Problem}
\newcommand{\R}{\mathbb{R}}
\newcommand{\C}{\mathbb{C}}
\newcommand{\tr}{\operatorname{tr}}

\newcommand{\dd}{\,\mathrm{d}}

\newcommand{\rank}{\operatorname{rank}}
\title{\textbf{A Power Geometry and Hypergeometric Functions}}
\author{Shintaro Yoshizawa\\[2mm]\small Nagoya Mathematical and Information Science Research\\\small\href{mailto:shintaro.yoshizawa.net@gmail.com}{shintaro.yoshizawa.net@gmail.com}}
\date{}
\begin{document}
\maketitle

\section{Introduction}
The present paper is a corrected English translation of the original Japanese source cited as \cite{original}. Corrections are limited to typographical, bibliographical, and explicitly identified mathematical notation issues; the structure and substantive mathematical content of the original are otherwise retained.

The term ``power geometry'' is not a conventional technical term. Rather, it expresses a research programme that I have pursued, centred on the \emph{power transformation}, \emph{deformation}, and \emph{duality} of functions (in the sense of optimization theory). My intention is to investigate pure mathematics within applied mathematics by means of analytic, geometric, and algebraic methods, and, where appropriate, to return the resulting theory to applications.

There are two principal motivations for studying power geometry. The first arose from a simple question that has occupied me since around 1996: ``How are probability distributions generated?'' --- the origin of probability distributions. Since information geometry \cite{3} has established the view that a family of probability distributions itself constitutes a space, the question may also be interpreted, somewhat grandiosely, as ``What is the origin of a space?''

Among the various ways of constructing families of probability distributions, I have been particularly interested in the following:
\begin{enumerate}[label=(\roman*)]
\item the maximum-likelihood principle together with a definition of the mean that attains it;
\item identifying parameters of rational functions and constructing finite measures; and
\item constructing distributions as solutions of ordinary differential equations in the frequency domain.
\end{enumerate}
The first was considered by Gauss (1809), Poincar\'e (1912), Keynes (1911), and others. Laws of error associated with the arithmetic, geometric, and harmonic means were discussed, and Matsunawa (1994) further developed Keynes's work from a modern viewpoint \cite{15}. In statistics, Pearson (1895) proposed the Pearson system, deriving families of probability distributions from an ordinary differential equation whose logarithmic-likelihood derivative is prescribed as a rational function. In mathematical control, the Pearson system was further generalized as the geometry of rational functions and developed into an exponent theorem for rational maps from the Riemann sphere to Grassmann manifolds \cite{8}. The third direction, I believe, may offer an especially interesting development from the viewpoint of analytic geometry.

Beginning with Meixner's derivation of families of orthogonal polynomials from generating functions \cite{16}, Laha and Lukacs characterized the Meixner class as probability distributions by identifying the coefficient parameters of a nonlinear ordinary differential equation modelling a conditional quadratic regression problem \cite{13}. My observations of the nonlinear ODE of Laha and Lukacs, together with my participation in the 1999 Hokkaido University workshop ``On the Schwarzian Derivative,'' led me to focus on the relation between Schwarzian equations and probability distribution functions.

The second motivation for power geometry arose from the problem of determining eigenvalues (or singular values), matrix inverses, and matrix factorizations by means of dynamical systems whose states are matrices, not necessarily square matrices. The fundamental example is the optimization problem considered by von Neumann (1937): for symmetric matrices $A$ and $B$, determine the extrema of
\[
f(X)=\tr(AXBX^T),\qquad X\in SO(n).
\]
Brockett \cite{6} derived a double-Lax dynamical system as the gradient flow of $f$ with respect to the Killing form on the special orthogonal group, and clarified its character as an analogue of an eigenvalue-computation flow for a symmetric matrix, as well as its interpretation as an eigenvalue-sorting flow \cite{6}. Taking the constant symmetric matrix $B$ to be diagonal with distinct eigenvalues, Brockett \cite{6,7} obtained the following gradient flow on $SO(n)$ for computing the eigenvalues of the constant symmetric matrix $A$:
\begin{equation}\label{eq:brockett}
\dot X=AXB-XBX^TAX,\qquad X\in SO(n).
\end{equation}
Here $X^T$ denotes the transpose of $X$. Upon setting $L=XBX^T$, one obtains the differential equation on symmetric matrices
\begin{equation}\label{eq:lax}
\dot L=[L,[L,A]],
\end{equation}
where $[Z_1,Z_2]=Z_1Z_2-Z_2Z_1$. For computation, \eqref{eq:lax} is easier to implement because its state space is linear. Nevertheless, if the symmetric matrices $A$ and $B$ are restricted to be positive definite, one can prove that \eqref{eq:brockett} remains a gradient flow even when its domain is extended to all rectangular matrices. Moreover, if $X$ is an $n\times k$ rectangular matrix, $k\le n$, one obtains a gradient flow that simultaneously computes the $k$ largest eigenvalues and eigenvectors of a positive-definite symmetric matrix $A$.

\begin{fact}[24]
Let $X\in\R^{n\times k}$, $k\le n$, and let $A\in\R^{n\times n}$ and $B\in\R^{k\times k}$ be positive-definite symmetric matrices. For tangent vectors $V_1,V_2\in T_X\R^{n\times k}$, define the Riemannian metric by
\[
\langle V_1,V_2\rangle=\tr(AV_1BV_2^T).
\]
With respect to this metric, the potential
\begin{equation}\label{eq:potential3}
f(X)=\frac14\tr\!\left((AXBX^T)^2\right)-\frac12\tr(A^2XB^2X^T)
\end{equation}
has negative gradient flow
\begin{equation}\label{eq:flow4}
\dot X=AXB-XBX^TAX,\qquad X\in\R^{n\times k}.
\end{equation}
Furthermore, when $B$ is diagonal, the flag manifold is both an attractor and an invariant manifold. If $B=I$, the system preserves the rank of $X$; hence, for a full-rank initial condition, $\rank X(0)=k$, the Stiefel manifold $\mathrm{St}(n,k)$ is invariant, i.e. $X(t)\in\mathrm{St}(n,k)$ for all $t\ge0$.
\end{fact}

In optimization theory, once a potential function can be written explicitly, it is natural to investigate its dual potential as well. When $B=I$, equation \eqref{eq:flow4} is a flow for finding the $k$-dimensional principal subspace of $A$, namely the invariant subspace spanned by eigenvectors associated with its $k$ largest eigenvalues. Around 2004 there was a conjecture (Problem 3.9 in \cite{5}) that the gradient flow of the Legendre dual of the potential in \eqref{eq:potential3} with $B=I$ should yield a dynamical system for finding the $k$-dimensional minor subspace, spanned by the eigenvectors associated with the $k$ smallest eigenvalues. This was answered affirmatively in \cite{26}. For a general positive-definite symmetric matrix $B$, the problem remains open; the matrix $B$ constitutes an obstruction when a multivalued function is to be made single-valued.

There are other potential functions for finding the $k$ largest eigenvalues and eigenvectors of a positive-definite symmetric matrix $A$. Interestingly, by taking powers of functions of matrix variables and analytically deforming (twisting) them, one can connect, by a real parameter $\alpha$, a potential for the largest eigenvalues and eigenvectors of $A$ with one for the smallest eigenvalues and eigenvectors. One example is
\begin{equation}\label{eq:galpha}
g_\alpha(X)=\frac12\tr(X^TAX)-\frac12\tr\left(\frac{(X^TX+B)^\alpha-I}{\alpha}\right),\qquad X\in\R^{n\times k},
\end{equation}
where $A$ and $B$ are positive-definite symmetric matrices. As $\alpha\to2$, this becomes a potential whose gradient flow finds the $k$ largest eigenvalues and eigenvectors of $A$; as $\alpha\to0$, it becomes a potential whose gradient flow finds the $k$ smallest eigenvalues and eigenvectors of $A$.

Both the potential \eqref{eq:potential3} and \eqref{eq:galpha}, in their forms associated with the largest eigenvalues and eigenvectors, possess a difference-of-convex (DC) structure. In contrast, the potential in \eqref{eq:galpha} associated with the smallest eigenvalues and eigenvectors does not have a DC structure. The gradient flow for the minor subspace involves a topological issue different from that of the principal-subspace flow \cite{14}. The study of function families including \eqref{eq:galpha} remains a subject for future work.

In information geometry \cite{3}, the potential functions associated with Gaussian distributions have a DC structure when expressed in terms of mean and variance parameters. For exponential families, viewing the potential in canonical ($e$-) coordinates or expectation ($m$-) coordinates allows the DC function to be convexified, bringing out the naturalness of the dual-connection structure. The geometry of the parameter space of matrix-valued Gaussian random variables whose entries take values in the reals, complexes, quaternions, or octonions, together with its power geometry, is another subject for future research; this investigation has only just begun \cite{23}.

Finally, I mention a setting in which the two motivations for power geometry are linked by the Schwarzian equation. If the initial value of $X$ in \eqref{eq:flow4} is chosen on a flag manifold, then the change of variables $L=XBX^T$ yields a matrix Riccati differential equation. By considering a matrix-valued Schwarzian derivative and the correspondence between matrix Schwarzian equations and matrix Riccati equations, I am working toward a deeper understanding of dynamical systems for eigenvalues and singular values, with the aim of providing design principles for new classical or quantum information-processing algorithms.

The remainder of the paper is organized as follows:
\begin{itemize}
\item \S2: Generation of Schwarzian equations;
\item \S3: Schwarz divergences (kernels and Schwarzian derivatives); and
\item \S4: Matrix dynamical systems and matrix Schwarzian equations.
\end{itemize}

The present research programme is still developing. I express my sincere gratitude to Professor Hiroshi Matsuzoe, the principal investigator of the research project ``New Developments in the Geometry of Statistical Manifolds,'' and to all the institutions and colleagues involved, for the opportunity to present this work.

\section{Generation of Schwarzian Equations}
As a deductive approach to the origin of probability distributions, we revisit the nonlinear ODE considered by Laha and Lukacs, discuss its relation to the Schwarzian derivative, and outline problems in power geometry involving pseudodifferential operators. As an experimental approach, we repeatedly apply power transformations to functions, thereby generating families that may serve as candidates for probability densities; from these families we then generate Schwarzian equations and explore concretely the relation between power functions and Schwarzian equations.

\subsection{The nonlinear ODE of Laha and Lukacs}
Let $F(x)$ be a probability distribution function on the real line, let $f(x)$ be its probability density, and let $\phi(t)$ be its characteristic function, i.e.,
\[
\phi(t)=\int_{-\infty}^{\infty}e^{\sqrt{-1}xt}\dd F(x)
=\int_{-\infty}^{\infty}e^{\sqrt{-1}xt}f(x)\dd x.
\]

\begin{theorem}[12]
The Meixner class is characterized by the nonlinear ODE
\begin{equation}\label{eq:laha}
(1-a)\left(\frac{\phi'}{\phi}\right)'-a\left(\frac{\phi'}{\phi}\right)^2
=\sqrt{-1}\,b\frac{\phi'}{\phi}-c,
\end{equation}
and by the discriminant $\Delta$ of its solutions. Here $a,b,c\in\R$, $\phi'=\dd\phi/\dd t$, and $\Delta=b^2-4ac$. The cases are:
\begin{enumerate}[label=(\roman*)]
\item $\Delta=0$:
  \begin{enumerate}[label=(\alph*)]
  \item if $a=b=0$, then $f$ is Gaussian;
  \item if $a\ne0$ and $b\ne0$, then $f$ is Gamma.
  \end{enumerate}
\item $0<\Delta$:
  \begin{enumerate}[label=(\alph*)]
  \item if $a=0$, $b<0$, and $0<c$, then $f$ is Poisson;
  \item if $a\ne0$, then, according to the sign of $c$, $f$ is binomial or negative binomial.
  \end{enumerate}
\item $\Delta<0$:
  \begin{enumerate}[label=(\alph*)]
  \item if $0<a<1$, then $f$ is a Meixner hypergeometric distribution;
  \item if $a=1$, then $f$ is Cauchy.
  \end{enumerate}
\end{enumerate}
\end{theorem}

\begin{remark}
If
\[
\varphi(t)=\int^t\phi(s)\dd s,
\]
then \eqref{eq:laha} becomes
\begin{equation}\label{eq:laha2}
(1-a)\left(\frac{\varphi''}{\varphi'}\right)'-a\left(\frac{\varphi''}{\varphi'}\right)^2
=\sqrt{-1}\,b\frac{\varphi''}{\varphi'}-c.
\end{equation}
For $a=1/3$, the left-hand side of \eqref{eq:laha2} is $(2/3)$ times the Schwarzian derivative
\[
\{\varphi;t\}=\left(\frac{\varphi''}{\varphi'}\right)'-\frac12\left(\frac{\varphi''}{\varphi'}\right)^2.
\]
With the change of variable $\psi=\phi'/\phi$, equation \eqref{eq:laha} is plainly a Riccati equation for $\psi$ and is therefore explicitly solvable.
\end{remark}

\begin{fact}[3]
If $p(x,t)=ax^2+bx+c$, then \eqref{eq:laha} is equivalent to
\begin{equation}\label{eq:bolger}
\phi(t)\frac{\dd^2}{\dd t^2}\log\phi(t)
=-\int_{-\infty}^{\infty}e^{\sqrt{-1}tx}p(x,t)f(x)\dd x.
\end{equation}
\end{fact}

\begin{problem}
The function $p(x,t)$ on the right-hand side of \eqref{eq:bolger} resembles an amplitude function in an oscillatory integral. Generalize $p(x,t)$ to a power function that contains quadratic functions of $x$ as a special case, and investigate the correspondence between the amplitude function and the dual geometric structure of probability distributions. The matrix-valued version, in which $x$ is replaced by a matrix variable, should also be studied.
\end{problem}

\subsection{Generation of Schwarzian equations by multiple power transformations}
When one wishes to analyze skewed data under an assumption of normality, the Box--Cox transformation is used to bring a random variable closer to a Gaussian distribution. Specifically, one estimates $\alpha$ so that the transformation
\[
L_\alpha(x)=\frac{1-x^\alpha}{\alpha}
\]
brings a skewed random variable $x$ closer to Gaussianity. More generally, define, for $x>0$,
\[
\log_q(x)=\frac{x^{1-q}-1}{1-q},
\qquad
\exp_q(x)=\{1+(1-q)x\}^{1/(1-q)},
\]
whenever the expression in braces is positive. Then
\[
\lim_{\alpha\to0}L_\alpha(x)=-\log x,
\qquad
\lim_{q\to1}\log_q(x)=\log x,
\qquad
\lim_{q\to1}\exp_q(x)=\exp x.
\]

We repeatedly twist a family of power functions by repeated $L_\alpha$ transformations, thereby constructing deformed multiple-power families. Computing the Schwarzian derivatives of the resulting families yields concrete Schwarzian equations. The following is an example (the arguments of all powers are assumed nonnegative):
\begin{equation}\label{eq:la}
L_a(x)=t\quad\Longrightarrow\quad
x=
\begin{cases}
(1-at)^{1/a},&a\ne0,\\
\exp(-t),&a=0.
\end{cases}
\end{equation}
Furthermore,
\begin{equation}\label{eq:lbla}
L_b(L_a(x))=t\quad\Longrightarrow\quad
x=
\begin{cases}
\{1-a(1-bt)^{1/b}\}^{1/a},&a,b\ne0,\\
\exp\{-(1-bt)^{1/b}\},&a=0,\ b\ne0,\\
\{1-a\exp(-t)\}^{1/a},&a\ne0,\ b=0,\\
\exp\{-\exp(-t)\},&a=b=0.
\end{cases}
\end{equation}

\begin{definition}
The deformed multiple-power family is defined by
\[
(L_{a_n}\circ\cdots\circ L_{a_2}\circ L_{a_1})(x)
=(L_{b_m}\circ\cdots\circ L_{b_2}\circ L_{b_1})(t),
\]
where $a_i,b_i\in\R$. Choosing the range of $t$ so that all powers have nonnegative arguments, and solving for $x$, we write $x\in F_{n,m}(t)$.
\end{definition}

These families form the following array:
\[
\begin{array}{cccccc}
\vdots&\vdots&\vdots&\\
F_{0,2}(t)&\longrightarrow&F_{1,2}(t)&\longrightarrow F_{2,2}(t)\longrightarrow\cdots\\
\uparrow&&\uparrow&&\uparrow\\
F_{0,1}(t)&\longrightarrow&F_{1,1}(t)&\longrightarrow F_{2,1}(t)\longrightarrow\cdots\\
\uparrow&&\uparrow&&\uparrow\\
F_{0,0}(t)&\longrightarrow&F_{1,0}(t)&\longrightarrow F_{2,0}(t)\longrightarrow\cdots
\end{array}
\]

\begin{definition}
For $x\in F_{n,m}(t)$, let $\{x;t\}$ denote an element of $S_{n,m}(t)$; that is, $\{x;t\}\in S_{n,m}(t)$.
\end{definition}

\begin{theorem}[H. Schwarz, 1872]
Consider the Gauss hypergeometric equation
\begin{equation}\label{eq:hypergeom}
\frac{\dd^2z}{\dd t^2}
+\frac{\gamma-(\alpha+\beta+1)t}{t(1-t)}\frac{\dd z}{\dd t}
-\frac{\alpha\beta}{t(1-t)}z=0,
\end{equation}
where $\alpha,\beta,\gamma\in\C$. If $z_1(t),z_2(t)$ are linearly independent solutions and $s(t)=z_1(t)/z_2(t)$, then the corresponding Schwarzian equation is
\begin{equation}\label{eq:schwarz-hyper}
\{s;t\}=\frac{1-\lambda^2}{2t^2}
+\frac{1-\mu^2}{2(1-t)^2}
+\frac{1+\nu^2-\lambda^2-\mu^2}{2t(1-t)},
\end{equation}
where
\[
\lambda=1-\gamma,\qquad
\mu=\gamma-\alpha-\beta,\qquad
\nu=\alpha-\beta.
\]
\end{theorem}

\begin{fact}
Let
\[
x=\{1-a(1-b(pt+q))^{1/b}\}^{1/a}\in F_{2,0}(pt+q),
\qquad p\ne0,
\]
with $p,q\in\R$. For suitable special choices of the parameters, this family belongs to the family of Schwarzian derivatives determined by the Gauss hypergeometric equation. In fact,
\begin{equation}\label{eq:schwarz-family}
\{x;t\}=
\frac{p^2\left[(a^2b^2-1)(-bpt-bq+1)^{2/b}
+2a(1-b^2)(-bpt-bq+1)^{1/b}+b^2-1\right]}
{2(bpt+bq-1)^2\{a(-bpt-bq+1)^{1/b}-1\}^2}.
\end{equation}
When $b=1$, the denominator is a polynomial of degree four in $t$, while the numerator is of degree two. This can be verified further by taking $ap=-1$ and $q=1$.
\end{fact}

\begin{remark}
If $a=\pm1$, $b=0$, $p=1$, and $q=0$, then $\{x;t\}=-1/2$.
\end{remark}

\begin{problem}
For the families $F_{n,m}$, consider the families of Schwarzian derivatives $S_{n,m}$ and investigate the structure of the collection $\{S_{n,m}\}_{n,m}$.
\end{problem}

\begin{problem}
Consider the matrix zeta function
\[
M_\alpha(X)=\frac{X^{-\alpha}-I}{\alpha},
\]
where $X$ is square and $\alpha\in\R$. Construct deformed multiple-power matrix zeta functions and investigate the corresponding Schwarzian equations. For operator zeta functions and matrix zeta functions, see, for example, \cite{9}.
\end{problem}

Finally, let us power-deform the general solution of a simple special Schwarzian equation and examine its Schwarzian derivative.

\begin{fact}
If
\[
\{x;t\}=-2k^2,
\]
where $k\in\R$ is constant, then the general solution is
\begin{equation}\label{eq:general-mobius}
x(t)=\frac{a\exp(kt)+b\exp(-kt)}{c\exp(kt)+d\exp(-kt)},
\qquad ad-bc\ne0.
\end{equation}
Power-deforming the exponential terms gives
\begin{equation}\label{eq:xtau}
x_\tau(t)=
\frac{a(1+\tau kt)^{1/\tau}+b(1-\tau kt)^{1/\tau}}
{c(1+\tau kt)^{1/\tau}+d(1-\tau kt)^{1/\tau}},
\qquad ad-bc\ne0,
\end{equation}
and its Schwarzian derivative is
\begin{equation}\label{eq:xtau-schwarz}
\{x_\tau;t\}=
\frac{2(\tau^2-1)k^2}{(\tau kt+1)^2(\tau kt-1)^2}.
\end{equation}
\end{fact}

\begin{remark}
The result \eqref{eq:xtau-schwarz} is intriguing in relation to the classical univalence criteria of Nehari \cite{17} and Hille \cite{10}, and to the Ahlfors--Weill criterion for quasiconformal extension \cite{2}.
\end{remark}

\begin{problem}
Equation \eqref{eq:general-mobius} was used in Kobayashi's proof \cite{11} that, on a complete Riemannian manifold whose Ricci curvature is proportional to the Riemannian metric and negative, the projectively invariant pseudodistance agrees, up to a constant factor, with the distance induced by the given Riemannian metric. What is the geometric meaning corresponding to the deformed Schwarzian equation \eqref{eq:xtau-schwarz}?\footnote{See the section on affine connections in \cite{12}.}
\end{problem}

\section{Schwarz Divergence: Kernels and the Schwarzian Derivative}
Matrix-valued Schwarzian derivatives have roots, as far as I know, in several streams of research: the discovery of the Virasoro algebra by Gelfand and Fuchs in 1967 (see, e.g., Chapter 7 of Ovsienko and Tabachnikov \cite{18}); work arising from geometric optimal control theory (Zelikin \cite{27} and the references to Agrachev and others in Paiva and Dur\'an \cite{19}); and the line of work in function theory associated with B. Schwarz \cite{21}.

In this section, rather than following these approaches directly, we revisit the Schwarzian derivative through its relation to a kernel function and derive a matrix-valued Schwarzian derivative in a natural way. The main idea is as follows. In information geometry, divergences play a fundamental role. We reinterpret the Bregman divergence generated by a convex function through the logarithmic mean rate of change of a convex monotone function $f$, and use this to define a Schwarz divergence. The Schwarz functional is defined on the Cartesian square of the domain of $f$ and characterized so that, on the diagonal, it becomes the Schwarzian derivative. Finally, a matrix-valued Schwarzian derivative is constructed from the matrix Schwarz functional.

\subsection{Schwarz Divergence}
Let $I$ be a nonempty open interval of the real line, and let $f$ be at least three times continuously differentiable and strictly convex on $I$. The Bregman divergence $D_f(x,y)$ on $I\times I$ is fundamental in convex optimization. The following facts are standard.

\begin{fact}
Define
\begin{equation}\label{eq:bregman}
D_f:I\times I\ni(x,y)\mapsto f(x)-f(y)-\partial_yf(y)(x-y)\in[0,\infty).
\end{equation}
Then:
\begin{enumerate}[label=(\roman*)]
\item $D_f(x,y)=f(x)+f^*(w)-xw$;
\item $0\le D_f(x,y)$;
\item $D_f(x,y)=0$ implies $x=y$,
\end{enumerate}
where $w=\partial_yf(y)$ and $f^*$ denotes the Legendre dual of $f$.
\end{fact}

Suppose that $f$ also satisfies the monotonicity condition. Define the Schwarz divergence by
\begin{equation}\label{eq:schwarz-div}
D_f^{\mathrm S}:I\times I\ni(x,y)\mapsto
\frac{1}{x-y}-\frac{\partial_yf(y)}{f(x)-f(y)}\in[0,\infty).
\end{equation}

\begin{fact}
For \eqref{eq:schwarz-div}:
\begin{enumerate}[label=(\roman*)]
\item $0\le D_f^{\mathrm S}(x,y)$;
\item $D_f^{\mathrm S}(x,y)=0$ implies $x=y$;
\item $\partial_xD_f^{\mathrm S}(x,y)=\partial_{xy}^2K_f(x,y)$ for $x\ne y$,
\end{enumerate}
where
\begin{equation}\label{eq:Kf}
K_f(x,y)=\log\frac{f(x)-f(y)}{x-y},\qquad x\ne y.
\end{equation}
We call $K_f(x,y)$ the logarithmic Loewner kernel and $\partial_{xy}^2K_f(x,y)$ the Schwarz functional.
\end{fact}

The following identity is well known:
\begin{fact}
\begin{equation}\label{eq:diag-schwarz}
\partial_{xy}^2K_f(x,y)
=\frac{f'(x)f'(y)}{(f(x)-f(y))^2}-\frac{1}{(x-y)^2}
\xrightarrow[y\to x]{}\frac16\{f;x\},
\end{equation}
where
\[
\{f;x\}=\left(\frac{f''}{f'}\right)'-\frac12\left(\frac{f''}{f'}\right)^2.
\]
\end{fact}
This identity motivates the terminology ``Schwarz divergence.''

To state the relation between Loewner matrices and the Schwarz functional, we first recall the definition of a Loewner matrix.

\begin{definition}
Let $f$ be defined on $(a,b)$ and let $x_1,\ldots,x_n\in(a,b)$ be distinct. Define
\begin{equation}\label{eq:loewner}
L^{(n)}(x_1,\ldots,x_n;f)_{ij}=
\begin{cases}
\dfrac{f(x_i)-f(x_j)}{x_i-x_j},&i\ne j,\\[2mm]
f'(x_i),&i=j.
\end{cases}
\end{equation}
\end{definition}

\begin{fact}
The Loewner matrix and the Schwarz functional are related by
\begin{equation}\label{eq:loewner-schwarz}
\partial_{xy}^2K_f(x,y)\,(f(x)-f(y))^2
=\det L^{(2)}(x,y;f).
\end{equation}
\end{fact}

\subsection{Matrix-valued Schwarzian Derivatives}
We begin with the definition of a matrix-monotone function in order to construct a matrix-valued Schwarzian derivative from a matrix Schwarz functional.

\begin{definition}
Let $f(t)$ be a real-valued continuous function on an open interval $I$. For a Hermitian matrix whose eigenvalues all lie in $I$, the matrix $f(A)$ is defined by functional calculus. If, for matrices $A,B$ of every size,
\[
A\le B\quad\Longrightarrow\quad f(A)\le f(B),
\]
then $f$ is called matrix monotone.
\end{definition}

\begin{remark}
It is known that if a real-valued continuous function on the whole real line is matrix monotone for all $2\times2$ Hermitian matrices, then $f$ is affine [21, p.~87].
\end{remark}

\begin{definition}
Let $X,Y$ be Hermitian matrices with spectra in $(a,b)$ and let $f$ be increasing and matrix monotone. Define the matrix-variable logarithmic Loewner kernel by
\begin{equation}\label{eq:matrix-kernel}
K_f(X,Y)=6\,\log\left\{(f(X)-f(Y))(X-Y)^{-1}\right\}.
\end{equation}
\end{definition}

Let $\{X(t)\}$ be a smooth curve. Consider the Taylor expansion of
\[
\partial_{ts}^2K_f(X(t),X(s))
=\frac{\partial^2K_f}{\partial t\,\partial s}(X(t),X(s))
\]
and take the matrix-valued limit as $t\to s$. We define the result as the matrix Schwarzian derivative $S(Z(t))$, or $\{Z;t\}$.\begin{equation}\label{eq:matrix-schwarz}
S(Z(t))=\left\{(Z'(t))^{-1}Z''(t)\right\}'
-\frac12\left\{(Z'(t))^{-1}Z''(t)\right\}^2,
\end{equation}
where $Z=Z(t)\in\R^{n\times n}$ and $\det Z(t)\ne0$ for $t\in I$.

\section{Matrix Dynamical Systems and Schwarzian Equations}
Matrix Schwarzian equations are related to matrix Riccati equations through Hamiltonian systems \cite{27,28}. In this section we investigate the relation between the principal-subspace flow
\begin{equation}\label{eq:oja}
\dot X=(I-XX^T)AX,\qquad X\in\R^{n\times k},\quad k\le n,
\end{equation}
for a positive-definite symmetric matrix $A$, and matrix Schwarzian equations, together with intrinsic properties of the matrix Schwarzian derivative itself.

\begin{proposition}[Zelikin, 27]
Let $A$ and $B$ be symmetric. If $Z(t)$ is a solution of the matrix Schwarzian equation
\begin{equation}\label{eq:zelikin-schwarz}
S(Z(t))=2(B(t)-A'(t)),
\end{equation}
then the matrix variable
\begin{equation}\label{eq:W}
W(t)=-\frac12(Z'(t))^{-1}Z''(t)-A(t)
\end{equation}
is a solution of the Riccati equation
\begin{equation}\label{eq:riccati}
\dot W=-B-AW-WA-W^2.
\end{equation}
Conversely, if $W$ solves \eqref{eq:riccati}, then any function $Z(t)$ satisfying \eqref{eq:W} is a solution of \eqref{eq:zelikin-schwarz}.
\end{proposition}

To relate the matrix Riccati equations obtained from the principal-subspace flow to the Riccati equation in Proposition~1, consider the two changes of variables
\[
N=-2A^{1/2}XX^TA^{1/2},
\qquad
L=2A^{1/2}XX^TA^{1/2}.
\]
Then \eqref{eq:oja} is transformed respectively into
\begin{align}
\dot N&=AN+NA+N^2,\label{eq:N}\\
\dot L&=AL+LA-L^2.\label{eq:L}
\end{align}

\begin{fact}
The following equivalences hold:
\begin{enumerate}[label=(\roman*)]
\item If $N=\frac12(Z')^{-1}Z''-A$, then \eqref{eq:N} is equivalent to $S(Z)=-2A^2$.
\item If $L=-\frac12(Z')^{-1}Z''+A$, then \eqref{eq:L} is equivalent to $S(Z)=2A^2$.
\end{enumerate}
\end{fact}

\paragraph{Proof.}
For (i),
\[
N'=-\frac12(Z')^{-1}Z''(Z')^{-1}Z''+\frac12(Z')^{-1}Z''',
\]
and
\[
N^2=\frac14(Z')^{-1}Z''(Z')^{-1}Z''
-\frac12(Z')^{-1}Z''A-\frac12A(Z')^{-1}Z''+A^2.
\]
Consequently,
\begin{align*}
N'-N^2
={}&-\frac34\{(Z')^{-1}Z''\}^2
+\frac12(Z')^{-1}Z'''
+\frac12(Z')^{-1}Z''A\\
&+\frac12A(Z')^{-1}Z''-A^2\\
={}&\frac12S(Z)+A^2+NA+AN.
\end{align*}
Thus (i) follows, and (ii) is proved in the same manner.

We next ask what becomes of the preceding results when $W$ is invertible and $W$ is reciprocally related to $(Z')^{-1}Z''$.

\begin{fact}
Assume $A$ may be time-dependent or constant, and define
\[
\det W\ne0,
\qquad
W(t)=\left\{-\frac12(Z')^{-1}Z''-A\right\}^{-1}.
\]
Then
\begin{align}
\dot W&=I+AW+WA
&&\Longleftrightarrow&& S(Z)=-2(A^2+A'),\label{eq:W26}\\
\dot W&=I+AW+WA+WA^2W
&&\Longleftrightarrow&& S(Z)=-2A',\label{eq:W27}\\
\dot W&=I+AW+WA+W(A^2+A')W
&&\Longleftrightarrow&& S(Z)=0.\label{eq:W28}
\end{align}
\end{fact}

\paragraph{Proof.}
Differentiating the inverse gives
\[
\dot W=-W\left\{-\frac12(Z')^{-1}Z''-A\right\}'W.
\]
Using
\[
\left\{-\frac12(Z')^{-1}Z''-A\right\}'
=\frac12(Z')^{-1}Z''(Z')^{-1}Z''
-\frac12(Z')^{-1}Z'''-A',
\]
and substituting the identities obtained from the definition of $W$, one obtains
\begin{equation}\label{eq:Wmaster}
\dot W=I+AW+WA+W\left(A^2+A'+\frac12S(Z)\right)W.
\end{equation}
Equations \eqref{eq:W26}--\eqref{eq:W28} follow immediately.

We next examine the case in which $W$ and $(Z')^{-1}Z''$ are related by a fractional transformation, analogous to a Cayley transform.

\begin{fact}
Assume
\[
\det\left\{\frac12(Z')^{-1}Z''+A\right\}\ne0,
\qquad \det A\ne0,
\]
and define
\begin{equation}\label{eq:cayley}
W(t)=-\left\{\frac12(Z')^{-1}Z''-A\right\}
\left\{\frac12(Z')^{-1}Z''+A\right\}^{-1}.
\end{equation}
Then
\begin{equation}\label{eq:cayley-derivative}
\dot W=-\frac14(I+W)S(Z)A^{-1}(I+W)
-\frac12(I-W)A(I-W)
+\frac12A'(I-W)A^{-1}(I+W).
\end{equation}
\end{fact}

Equation \eqref{eq:cayley-derivative} yields
\begin{align}
\dot W&=-A-WAW+\frac12A'(I-W)A^{-1}(I+W)
&&\Longleftrightarrow&&S(Z)=2A^2,\\
\dot W&=AW+WA+\frac12A'(I-W)A^{-1}(I+W)
&&\Longleftrightarrow&&S(Z)=-2A^2,\\
\dot W&=-\frac12(I-W)A(I-W)+\frac12A'(I-W)A^{-1}(I+W)
&&\Longleftrightarrow&&S(Z)=0.
\end{align}
The proof is omitted.

Finally, we consider the matrix Schwarzian derivative from the viewpoint of Fanning frames.

\begin{definition}[18]
A smooth curve $c(t)$ of $n$-dimensional subspaces of $\R^{2n}$ is called \emph{Fanning} if, at every $t$, its tangent vector $\dot c(t)$ is an invertible linear map from $c(t)$ to the quotient space $\R^{2n}/c(t)$. A Fanning curve can be represented by a frame. If $A(t)$ is a smooth curve of rank-$n$ matrices of size $2n\times n$, so that $A$ is a frame, then the curve of $n$-dimensional subspaces generated by the columns of $A$ is Fanning if and only if the $2n\times2n$ matrix $(A,\dot A)$ is invertible for every $t$. Accordingly, a matrix $A$ for which $(A,\dot A)$ is invertible is called a Fanning frame.
\end{definition}

Paiva and Dur\'an \cite{19} characterize the Schwarzian derivative of a Fanning frame in terms of the coefficient of a second-order matrix differential equation (Theorem 3.4 of \cite{19}). For a Fanning frame written in the form used here, their matrix Schwarzian derivative agrees with
\[
\left\{(Z')^{-1}Z''\right\}'-\frac12\left\{(Z')^{-1}Z''\right\}^2.
\]

The following theorem characterizes the relation among the Fanning frame $A$, the fundamental endomorphism $F$, and the matrix Schwarzian derivative $\{A;t\}$. Here the Fanning frame is defined by the matrix Schwarzian derivative, and $F$ is regarded as a linear transformation of $\R^{2n}$ induced by the Fanning frame.

\begin{theorem}
Let $F(t)$ be the fundamental endomorphism associated with the Fanning frame $A$ by
\begin{equation}\label{eq:F}
F(t)=(A(t),\dot A(t))
\begin{pmatrix}0&I\\0&0\end{pmatrix}
(A(t),\dot A(t))^{-1}.
\end{equation}
Here
\[
\dot A=\frac{\dd A}{\dd t},
\qquad
A(t)=\begin{pmatrix}I\\Z\end{pmatrix},
\qquad
\{A;t\}=\frac{\dd}{\dd t}\bigl((Z')^{-1}Z''\bigr)
-\frac12\bigl((Z')^{-1}Z''\bigr)^2.
\]
Then
\begin{equation}\label{eq:Fidentity}
\left(\frac{\dd F(t)}{\dd t}\right)^2A(t)=2A(t)\{A(t);t\}.
\end{equation}
\end{theorem}
The proof is omitted.

\begin{problem}
Generalize the matrix
\[
\begin{pmatrix}0&I\\0&0\end{pmatrix}
\]
defining the fundamental endomorphism $F$ to a power matrix, and generalize the identity matrix defining the Fanning frame to a symmetric matrix $A$ or related structures. Determine what properties the analogue of \eqref{eq:Fidentity} possesses for the resulting matrix $A$.
\end{problem}

\end{document}